\documentclass[11pt,a4paper]{article}
\usepackage{amssymb}
\usepackage{amsmath}
\newcommand{\ep}{p_{\epsilon}}
\newcommand{\eq}{q_{\epsilon}}
\newcommand{\et}{t_{\epsilon}}

\newcommand{\lcg}{L}
\newcommand{\lloop}{l}
\newcommand{\inj}{{\rm inj}}

\newcommand{\qed}{\hfill $\Box$}
\newcommand{\qeda}{\hfill $\Box$\\[1ex]}

\newcommand{\dll}{\pi (1+\lambda^{-1})}

\newcommand{\Cut}{{\rm Cut}}
\newcommand{\R}{\mathbb{R}}

\newtheorem{lemma}{Lemma}
\newtheorem{remark}{Remark}
\newtheorem{proposition}{Proposition}
\usepackage{amssymb}
\author{Hans-Bert Rademacher}
\title{The length of a shortest geodesic loop}
\date{}
\begin{document}
\maketitle
\begin{abstract}
We give a lower bound for the 
length of a non-trivial geodesic loop
on a simply-connected
and compact manifold 
of even dimension with a non-reversible Finsler metric 
of positive flag curvature. Harris and
Paternain use this estimate 
in their recent paper~\cite{HP} to give a geometric 
characterization of dynamically convex Finsler metrics
on the $2$-sphere.
\\[1ex]
\end{abstract}
On  compact and simply-connected Riemannian manifold
with positive sectional curvature $0<K\le 1$
the length of a non-constant geodesic loop is
bounded from below by $2 \pi.$ This result is due
to Klingenberg \cite{Kl} and is of importance in
proofs of the classical sphere theorem.

For a compact manifold $M$ with non-reversible Finsler 
metric $F$ the author introduced in 
\cite{Ra1} the {\em reversibility}
$\lambda:=\max \{F(-X); F(X)=1\}\ge 1.$ 
In this short note we show how one can use 
the results and methods from
\cite{Ra1} to obtain the following estimate for the 
length of a
geodesic loop depending on the flag curvature and the
reversibility.
\begin{proposition}
\label{thm:main}
Let $M$ be a compact and simply-connected differentiable 
manifold of even dimension $n \ge 2$
equipped with a non-reversible Finsler metric $F$ 
and flag curvature $K$ satisfying
$0<K\le 1.$ Then 
the length $l$ of a shortest non-constant
geodesic loop is bounded from below:
$l \ge \pi\left(1+\lambda^{-1}\right)\,.$
In addition the injectivity radius satisfies:
$\inj \ge \pi /\lambda.$
\end{proposition}
In \cite[Theorem 4]{Ra1} it is shown that
with the same assumptions
the length of a closed geodesic $c$
satisfies this estimate. Therefore Proposition~\ref{thm:main}
follows from Proposition~\ref{pro:lcg} which we
are going to prove in this note. Proposition~\ref{thm:main}
answers 
a question posed to the author by G. Paternain.
Using results by Hofer, 
Wysocki and Zehnder \cite{HWZ} and the statement of
Proposition~\ref{thm:main}
Harris and Paternain obtain the following 
geometric characterization of 
dynamically convex Finsler metrics on the $2$-sphere:
\begin{proposition}{\rm (Harris-Paternain \cite[Theorem B, Section 6]{HP})}\\
\label{thm:two-or-inf}
Let $F$ be a non-reversible Finsler metric on the 
$2$-sphere with reversibility $\lambda$ and flag
curvature $$\left(1-\frac{1}{1+\lambda}\right)^2<K\le 1.$$ 
Then the Finsler metric
is {\em dynamically convex,} in 
particular there are either two
geometrically distinct closed geodesics or there are
infinitely many geometrically distinct ones.
\end{proposition}
For existence results for closed geodesic of 
Finsler metrics on the $2$-sphere we refer to the 
recent survey \cite{Lo}
by Long and to \cite{Ra3}.
On the $n$-sphere $S^n$
there is a $1$-parameter family 
$F_{\epsilon}, \epsilon \in [0,1)$
of Finsler metrics (called
{\em Katok metrics}) with the
following properties: $F_0$ is the standard Riemannian
metric, for every $\epsilon \in (0,1)$ the 
metric is a non-reversible Finsler metric of 
constant flag curvature $1,$ the reversibility
is $\lambda=(1+\epsilon)/(1-\epsilon)$ and
the shortest geodesic loop is a closed geodesic
of length $\pi (1+\lambda^{-1}).$ This shows that the
estimate given in Proposition~\ref{thm:main} is 
sharp. In addition the number of closed geodesics
for $n=2$ is two if $\epsilon$ is irrational.
It is an open problem whether there is a non-reversible
Finsler metric on $S^2$ with a finite number
$N>2$ of geometrically distinct closed geodesics.
 \\[2ex]
We use the following notations on
a compact manifold $M$ with Finsler metric $F$ 
introduced in \cite{Ra1}:
For points $p,q \in M$ let
$\theta (p,q)$ be the minimal length of a 
piecewise differentiable curve $c:[0,1]\rightarrow M$
joining $p=c(0)$ and $q=c(1).$ For a non-reversible Finsler
metric we have in general $\theta(p,q)\not=\theta(q,p),$
i.e. $\theta$ defines in general a non-symmetric metric
on $M.$
Therefore $\theta (p,q)$ equals the length 
of a minimal geodesic
$c:[0,1]\rightarrow M$ joining $p=c(0)$ and 
$q=c(1)$, i.e. a geodesic with $L(c)=\theta(c(0),c(1)).$
Then $d:M\times M\rightarrow \R; d(p,q)=\left(\theta(p,q)+\theta(q,p)\right)/2$
defines a symmetric metric on $M.$ For a point $p \in M$
and an unit vector $v \in T_pM; F(v)=1$ let 
$c_v:\R \rightarrow M$
be the geodesic with $p=c(0); v=c'(0)$ and define
$t_v:=\sup \{t >0; \theta(c(0),c(t))=t\}.$ Then
$c_v(t_v)$ is a cut point of the point $p$
and the cut locus
$\Cut (p)$ of the point $p$ is given by
$\Cut(p)=\{c_v(t_v); v \in T_pM, F(v)=1\}.$
\\
We consider the following invariants and their
relations: 
The {\em symmetrized injectivity radius}
$d:=\inf \{d(p,q)\,;\, q \in {\rm Cut}(p)\},$
the length
$L$ of a shortest (nontrivial) closed geodesic and the length
$l$ of a shortest (nontrivial) geodesic loop.
\begin{lemma}
\label{lem:ld}
Let $(M,F)$ be a compact Finsler manifold, then the symmetrized injectivity
radius $d$ and the length $l$ resp. $L$ of a shortest non-trivial geodesic loop
resp. closed geodesic satisfy:
$2 d\le l \le L.$
\end{lemma}
{\sc Proof.}
Let $c:[0,l]\rightarrow M$ be a shortest geodesic loop
parametrized by arc length with $p=c(0)=c(l).$ Let $q=c(t), t \in (0,l)$
be the cut point, i.e. $c|[0,t]$ is minimal. It follows that
$l=L(c)\ge 
\theta(p,q)+\theta(q,p) = 2\,d(p,q) \ge 2 d.$ The inequality $L\ge l$ is obvious.
\qeda
The next ingredient in the Proof of Proposition~\ref{pro:lcg} is the following result:
\begin{lemma}
\label{lem:eins} {\rm \cite[Lemma7]{Ra1}}
Let $(M,F)$ be a compact Finsler manifold with reversibility 
$\lambda,$ symmetrized injectivity radius
$d$ and flag curvature $K\le1.$ 
If $2 d < \dll$ then 
the length 
$\lloop$ of a shortest 
non-trivial geodesic loop satisfies: $\lloop= 2\,d.$
\end{lemma}
With the help of these two Lemmata we prove the following
\begin{proposition}
\label{pro:lcg}
Let $(M,F)$ be a compact manifold with Finsler metric $F$ with 
reversibility $\lambda$ and flag curvature $K\le 1.$
If the symmetrized injectivity radius $d $ satisfies
$2 d<\dll$ then every shortest geodesic loop is
a closed geodesic, hence $\lcg=\lloop=2d.$
\end{proposition}
{\sc Proof.} Let $c:[0,l]\rightarrow M$ be a shortest geodesic loop parametrized
by arc length with $c(0)=c(l)=p.$ 
Let $q=c(t), t\in (0,l)$ be the cut point, i.e. $c|[0,t]$ is minimal.
We assume that $l < \dll.$ 
By Lemma~\ref{lem:ld} we obtain $2 d < \dll.$ 
Then we conclude from Lemma~\ref{lem:eins} that $\lloop=L(c)=2\,d.$
Since 
\begin{equation}
\label{eq:dd}
2 d = l =L(c)\ge
\theta(p,q)+\theta(q,p)= 2 d(p,q)
\end{equation}
and $q \in \Cut (p)$
it follows from
the definition of the symmetrized injectivity radius $d$ that
equality holds in Inequality~\ref{eq:dd}. Therefore $c\left|[t,l]\right.$ is a minimal
geodesic joining $q$ and $p.$
\\
For sufficiently small $\epsilon>0$ with $\ep=c(\epsilon)$
there is $\et\in (t,1)$ such that
$\eq=c\left(\et\right) \in \Cut(\ep),$
i.e. the geodesic $c|[\epsilon,\et] $ is minimal.
We conclude from the triangle inequality:
\begin{eqnarray*}
2\, d\left(\ep,\eq\right)=
\theta\left(\ep,\eq\right)+
\theta\left(\eq,\ep\right)
\le
\theta\left(\ep,\eq\right)+
\theta\left(\eq,p\right)+
\theta\left(p,\ep\right) =\\
=
\theta\left(p,q\right)+
\theta\left(q,p\right)
=L(c)=2\,d(p,q)\,.
\end{eqnarray*}
From the definition of the symmetrized injectivity radius
$d$ it follows that actually equality holds, 
i.e. the geodesic loop
is a closed geodesic.
\qeda
{\sc Proof of Proposition~\ref{thm:main}.}
We assume that $l < \dll\,$ and conclude from
Lemma~\ref{lem:ld}: $2 d=l<\dll.$
Then Proposition~\ref{pro:lcg} implies that $L=l=2 d < \dll.$
But in \cite[Theorem 4]{Ra1} it is shown that under 
the assumptions
of Proposition~\ref{thm:main} the length $L$ of a
shortest closed geodesic satisfies: $L\ge \dll.$ 
Therefore we obtain a contradiction, i.e. $l\ge \dll.$
\qed
\begin{remark}\rm
Under the assumptions of Proposition~\ref{pro:lcg} 
we have shown that
for any point $p \in M$ with a cut point $q \in M$ 
satisfying
$d(p,q)=d$ there is a shortest closed geodesic 
$c:[0,2d]\rightarrow M$
parametrized by arc length passing through $p$ and $q$, i.e.
$p=c(0)=c(2d); q=c(t); t \in (0,2d).$ Hence the restrictions
$c_1=c|\,[0,t]$ and $c_3=c|\,[t,2d]$ are minimal geodesics. The
cut point $q=c(t)$ is not a conjugate point 
since $t=\theta (p,q)<\pi$ and
$K \le 1.$ This implies that there is another 
minimal geodesic $c_2:[0,t]\rightarrow M$
joining $p$ and $q.$ Therefore Proposition~\ref{pro:lcg}
excludes the second case 
discussed in \cite[Remark 1]{Ra1} resp.
\cite[Lemma 9.7(b)]{Ra2}.
\end{remark}
\small

\small
{\sc Universit\"at Leipzig, Mathematisches Institut,
04081 Leipzig, Germany}\\[1ex]
E-Mail: {\tt rademacher@math.uni-leipzig.de}\\
{\tt http://www.math.uni-leipzig.de/\symbol{126}rademacher}
\end{document}